\magnification=\magstephalf
\input eplain
\input BMmacs
\input epsf

\centerline{{\titlefont SOME REMARKS ON MODIFIED DIAGONALS}}
\vskip .5cm

\centerline{{\byfont by}\quad {\namefont Ben Moonen and Qizheng Yin}}
\vskip 1cm

{\eightpoint 
\noindent
{\bf Abstract.} We prove a number of basic vanishing results for modified diagonal classes. We also obtain some sharp results for modified diagonals of curves and abelian varieties, and we prove a conjecture of O'Grady about modified diagonals on double covers.
\medskip

\par}

\section{Introduction}{Intro}

\ssection
Given a smooth projective variety~$X$ and a base point~$a$, Gross and Schoen introduced in~[\ref{GrSch}] modified diagonal cycles $\Gamma^n(X,a)$ on~$X^n$. For instance, $\Gamma^2(X,a) = \Delta_X - \bigl[X\times \{a\}\bigr] - \bigl[\{a\}\times X\bigr]$. In general, if $J \subset \{1,\ldots,n\}$ we define a closed subvariety $X^J \subset X^n$ by the condition that $x_i = a$ for all $i \notin J$; the modified diagonal $\Gamma^n(X,a)$ is then an alternating sum of the small diagonals on the~$X^J$. 

Gross and Schoen proved some vanishing results for the modified diagonals of curves, both in the Chow ring and modulo algebraic equivalence. In~[\ref{BeauVoi}], Beauville and Voisin prove that for a K3 surface~$X$ there is a distinguished point class $o_X \in \CH_0(X)$ and that $\Gamma^3(X,o_X) = 0$ in $\CH_2(X^3)$.  (Throughout we work with $\mQ$-coefficients.) A consequence of this is that the intersection pairing $\Pic(X)^{\otimes 2} \to \CH_0(X)$ takes values in $\mQ \cdot o_X$ and that $c_2(X) = 24\cdot o_X$.

Our interest in modified diagonals was sparked by the preprint~[\ref{OGr}] of O'Grady and the questions he asked. We were quickly able to answer one of these questions in the positive, proving that for a $g$-dimensional abelian variety~$X$ we have $\Gamma^m(X,a) = 0$ in $\CH_g(X^m)$ for $m>2g$ and any choice of base point; see~[\ref{MY}]. This result is included in the present note as Theorem~\ref{AVResult}. 

\ssection
In this paper, we give a simple motivic description of modified diagonals, and we collect a number of basic results about them. We also introduce and study some more general classes $\gamma^n_{X,a}(\alpha)$, for any $\alpha \in \CH(X)$, which for $\alpha = [X]$ give back the modified diagonals $\Gamma^n(X,a)$. 

We work over an arbitrary field and consider algebraic cycles (with $\mQ$-coefficients) modulo an adequate equivalence relation~$\sim$. We prove that $\Gamma^n(X,a) \sim 0$ if and only if the map $\gamma^n_{X,a}$ is zero modulo~$\sim$, and that this implies the vanishing of $\gamma^{n-s}_{X,a}(\alpha)$ for all classes~$\alpha$ in the image of the product map $\CH^{>0}(X)^{\otimes s} \to \CH(X)$. We also prove that if $f\colon X \to Y$ is surjective with generic fiber of dimension~$r$ then $\Gamma^n(X,a) \sim 0$ implies that $\Gamma^{n-r}\bigl(Y,f(a)\bigr) \sim 0$. Further we have a precise result about what happens when we change the base point (working in the Chow ring): if $\Gamma^n(X,a) = 0$ for some point~$a$ then for any other base point~$a^\prime$ we have $\Gamma^{2n-2}(X,a^\prime) = 0$.

In Section~\ref{CAV} we prove some sharp (or conjecturally sharp) vanishing results for modified diagonals of curves and abelian varieties. For a curve~$C$ we use the base point~$a$ to embed $C$ in its Jacobian~$J$. The vanishing of $\Gamma^n(C,a)$ is then equivalent to the vanishing of some components of the class $[C] \in \CH_1(J)$ with respect to the Beauville decomposition of $\CH(J)$. This is a problem that has been studied independently of modified diagonals, notably by Polishchuk, Voisin and the second author. Our Theorem~\ref{AVResult} about modified diagonals of abelian varieties, which proves a conjecture of O'Grady in~[\ref{OGr}], is an easy application of the results by Deninger and Murre~[\ref{DenMur}] about motivic decompositions of abelian varieties.

Finally, in Section~\ref{Double} we prove a conjecture of O'Grady about modified diagonals on double covers. We recently learned that Claire Voisin has proven a generalization of this result to covers of higher degree; this result is to appear in her forthcoming paper~[\ref{Voisin}].

\section{Definition and some basic properties of modified diagonals}{Basic}
\bigskip

\noindent
Throughout, Chow groups are taken with $\mQ$-coefficients.

\ssection
\ssectlabel{Motk}
Let $k$ be a field. Let $X$ and~$Y$ be smooth projective $k$-varieties. If $X$ is connected, let $\Corr_i(X,Y) = \CH_{\dim(X)+i}(X\times Y)$. In general, write~$X$ as a disjoint union of connected varieties, say $X = \coprod_\alpha\, X_\alpha$; then we let $\Corr_i(X,Y) = \oplus_\alpha\, \Corr_i(X_\alpha,Y)$. The elements of $\Corr_i(X,Y)$ are called correspondences from~$X$ to~$Y$ of degree~$i$. If $Z$ is a third smooth projective $k$-variety, composition of correspondences
$$
\Corr_i(X,Y) \times \Corr_j(Y,Z) \to \Corr_{i+j}(X,Z)
$$
is defined in the usual way: $(\theta,\xi) \mapsto \pr_{XZ,*}\bigl(\pr_{XY}^*(\theta) \cdot \pr_{YZ}^*(\xi)\bigr)$.

We denote by $\Mot_k$ the category of (covariant) Chow motives over~$k$. The objects are triples $(X,p,m)$ with $X$ a smooth projective $k$-variety, $p$ an idempotent in~$\Corr_0(X,X)$, and $m \in \mZ$. The morphisms from $(X,p,m)$ to $(Y,q,n)$ are the elements of
$$
q \circ \Corr_{m-n}(X,Y) \circ p
$$
(which is a subspace of $\Corr_{m-n}(X,Y)$), and composition of morphisms is given by composition of correspondences. The identity morphism on an object $(X,p,m)$ is $p \circ [\Delta_X] \circ p$, with $\Delta_X \subset X\times X$ the diagonal.

We have a covariant functor $h \colon \SmProj_k \to \Mot_k$, sending $X$ to $h(X) = (X,\Delta_X,0)$ and sending a morphism $f \colon X\to Y$ to the class of its graph $[\Gamma_f] \in \Corr_0(X,Y) = \Hom\bigl(h(X),h(Y)\bigr)$. We usually write $f_*$ instead of~$[\Gamma_f]$. 

There is a tensor product in~$\Mot_k$, making it into a $\mQ$-linear tensor category, such that $h(X) \otimes h(Y) = h(X\times Y)$. The unit object for this tensor product is the motive $\unitmot = h\bigl(\Spec(k)\bigr)$ of a point. If $M = (X,p,m)$ is an object of~$\Mot_k$ and $n\in \mZ$, we let $M(n) = (X,p,m+n)$. Then $M(n) = M \otimes \unitmot(n)$, and $\unitmot(1)$ is the Tate motive.

The Chow groups of a motive~$M$ are defined by $\CH_i(M) = \Hom_{\Mot_k}\bigl(\unitmot(i),M\bigr)$.

\ssection
\ssectlabel{h0h+Def}
Let $X$ be a connected smooth projective $k$-variety of dimension~$d$ with a rational point $a \in X(k)$. Then
$$
\pi_0 = X \times \{a\}
\quad\hbox{and}\quad
\pi_+ = [\Delta_X] - X \times \{a\}
$$
are orthogonal projectors, defining a decomposition
$$
h(X) = h_0(X) \oplus h_+(X)\, . \eqlabel{eq:h0h+}
$$
If there is a need to specify the base point, we use the notation $h_0(X,a)$ and $h_+(X,a)$.

If $f\colon X \to \Spec(k)$ is the structural morphism, $a \circ f$ is an idempotent endomorphism of~$X$ and $\pi_0$ is just the induced endomorphism $(a\circ f)_*$ of~$h(X)$. In particular, $f_* \colon h_0(X) \to h\bigl(\Spec(k)\bigr) = \unitmot$ is an isomorphism with inverse~$a_*$. On Chow groups we have $\CH\bigl(h_0(X)\bigr) = \mQ \cdot [a] \subset \CH(X)$.

\ssection
\ssectlabel{hXnDec}
We have a K\"unneth decomposition 
$$
h(X^n) = \bigl[h_0(X) \oplus h_+(X)\bigr]^{\otimes n} = \bigoplus_{J\subset \{1,\ldots,n\}}\; h_J(X^n)\, , \eqlabel{eq:hXnDec}
$$
where, for $J \subset \{1,\ldots,n\}$, we define
$$
h_J(X^n) =  h_{\nu_1}(X) \otimes \cdots \otimes h_{\nu_n}(X) \qquad \hbox{with}\quad \nu_i = \cases{+ & if $i \in J$\cr 0 & if $i \notin J$.} 
$$
The summand $h_{\{1,\ldots,n\}}(X^n) = h_+(X)^{\otimes n}$ will play a special role in what follows. Identifying $X^n \times X^n$ with $(X\times X)^n$, the projector onto this summand is~$\pi_+^{\otimes n}$.

\ssection
{\it Definition.\/} --- Retain the assumptions and notation of~\ref{h0h+Def}. For $n \geq 1$ define
$$
\gamma^n_{X,a} \colon h(X) \to h(X^n)
$$ 
by $\gamma^n_{X,a} = \pi_+^{\otimes n} \circ \Delta_{X,*}^{(n)}$, where $\Delta_X^{(n)} \colon X \to X^n$ is the diagonal morphism. We use the same notation $\gamma^n_{X,a}$ for the induced map on Chow groups $\CH(X) \to \CH(X^n)$ or on Chow groups modulo an adequate equivalence relation. Finally, we define
$$
\Gamma^n(X,a) \in \CH_d(X^n)
$$
(with $d=\dim(X)$) to be the image under $\gamma^n_{X,a}$ of the fundamental class $[X] \in \CH_d(X)$.

\ssection
\ssectlabel{GrossSchoen}
If $J$ is a subset of $\{1,\ldots,n\}$, we identify $X^J$ with the closed subvariety of~$X^n$ given by
$$
\bigl\{(x_1,\ldots,x_n) \in X^n \bigm| x_i = a\quad \hbox{if $i \notin J$}\bigr\}\, .
$$
Let $\phi_J = \phi_{X,J} \colon X^J \hookrightarrow X^n$ be the corresponding closed embedding. Let $\Delta^{(J)}_X \subset X^J$ be the small diagonal of~$X^J$, viewed as a cycle on~$X^n$.

If $\dim(X) = 0$ then $\gamma^n_{X,a}$ is the zero map. If $d = \dim(X)$ is positive, the cycle $\Gamma^n(X,a)$ is the modified diagonal cycle introduced by Gross and Schoen in~[\ref{GrSch}]. Explicitly, for $d>0$,
$$
\Gamma^n(X,a) = \sum_{\emptyset \neq J \subset \{1,\ldots,n\}}\, (-1)^{n-|J|}\cdot  \bigl[\Delta^{(J)}_X\bigr]\, .
$$

\ssection
\ssectlabel{h0h*h2d}
{\it Remark.\/} --- If $\dim(X) > 0$ we can refine~\eqref{eq:h0h+} to a decomposition
$$
h(X) = h_{2d}(X) \oplus h_\star(X) \oplus h_0(X)
$$
where $h_{2d}(X)$ and $h_\star(X)$ are the submotives of~$h(X)$ defined by the projectors $\pi_{2d} =  \{a\} \times X$ and $\pi_\star = [\Delta_X] - X\times\{a\} - \{a\} \times X$, respectively. For the study of modified diagonals this does not lead to a refinement, however, as for $n\geq 2$ the morphism $\gamma^n_{X,a} = \pi_+^{\otimes n} \circ \Delta_{X,*}^{(n)}$ is the same as the morphism $\pi_\star^{\otimes n} \circ \Delta_{X,*}^{(n)}$. To see this we have to show that 
$$
\bigl(\Delta_X^{(n)} \times \id_{X^n}\bigr)^* \pi_+^{\otimes n} = \bigl(\Delta_X^{(n)} \times \id_{X^n}\bigr)^* \pi_\star^{\otimes n}
$$ 
in $\CH(X \times X^n)$. (Use [\ref{DenMur}], Proposition~1.2.1.) Abbreviating $\Delta_X^{(n)}$ to~$\Delta$ and writing $p_i \colon X^n \to X$ for the $i$th projection, the difference $\bigl(\Delta_X^{(n)} \times \id_{X^n}\bigr)^* \bigl[\pi_+^{\otimes n} - \pi_\star^{\otimes n}\bigr]$ is a linear combination of terms 
$$
(\Delta \times \id_{X^n})^* (\beta_1 \otimes \cdots \otimes \beta_n) =  (\id_X \times p_1)^*\beta_1 \cdots (\id_X \times p_n)^*\beta_n
$$
where $\beta_1,\ldots,\beta_n \in \bigl\{\pi_{2d},\pi_\star\bigr\}$ and at least one~$\beta_j$ equals~$\pi_{2d}$. Now note that
$$
(\id_X \times p_i)^* \pi_{2d} \cdot (\id_X \times p_j)^* \pi_{2d} = 0\, ,
\quad\hbox{and}\qquad
(\id_X \times p_i)^* \pi_{2d} \cdot (\id_X \times p_j)^* \pi_\star = 0
$$ 
for all $i\neq j$.

\ssection
\ssectlabel{f_*Prop}
{\it Proposition. --- Let $f \colon X \to Y$ be a morphism of connected smooth projective $k$-varieties. Let $a \in X(k)$ and let $b = f(a)$. 

{\rm (\romno1)} The morphism $f_* \colon h(X) \to h(Y)$ is the direct sum of two morphisms $h_0(X,a) \to h_0(Y,b)$ and $h_+(X,a) \to h_+(Y,b)$.

{\rm (\romno2)} We have $\gamma^n_{Y,b} \circ f_* = f^{\otimes n}_* \circ \gamma^n_{X,a}$ for all $n \geq 1$.

{\rm (\romno3)} Suppose $f$ is generically finite of degree~$N$. Then $N \cdot \Gamma^n_{Y,b} = f^{\otimes n}_* \bigl(\Gamma^n_{X,a}\bigr)$ for all $n \geq 1$. 
\par}
\medskip

\Proof
For (\romno1), if $g \colon Y \to\Spec(k)$ is the structural morphism then $\pi_0(Y,b) = b_* \circ g_* = f_* \circ a_* \circ g_*$ and $\pi_0(X,a) = a_* \circ g_* \circ f_*$. Hence $\pi_0(Y,b) \circ f_* = f_* \circ \pi_0(X,a)$, and because $\pi_+ = \id - \pi_0$ also $\pi_+(Y,b) \circ f_* = f_* \circ \pi_+(X,a)$. Part~(\romno2) readily follows from this and (\romno3) follows by applying~(\romno2) to the class~$[X]$.
\QED

\section{Some vanishing results}{Firstvanish}

\ssection
\ssectlabel{modaer}
In what follows, we consider an adequate equivalence relation~$\aer$ on algebraic cycles, as in [\ref{YAIntro}], Section~3.1, and  we write $\Mot_{k,\aer}$ for the corresponding category of motives. If $M$ is an object of~$\Mot_{k,\aer}$, let $\A_i(M) = \Hom_{\Mot_{k,\aer}}\bigl(\unitmot(i),M\bigr)$ and $\A(M) = \bigoplus_{i\in\mZ}\, \A_i(M)$. In particular, if $X$ is a smooth projective $k$-variety, $\A_i(X) = \CH_i(X)/\aer$.

Given a connected smooth projective $k$-variety~$X$ with base point $a\in X(k)$, the decomposition~\eqref{eq:hXnDec} induces a decomposition
$$
\A(X^n) =  \bigoplus_{J \subset \{1,\ldots,n\}}\; \A_J(X^n)\, .
$$
This decomposition in general depends on the chosen base point.

Define a grading $\A(X^n) = \A_{\age{0}}(X^n) \oplus \cdots \oplus \A_{\age{n}}(X^n)$ by letting $\A_{\age{m}}(X^n)$ be the sum of all $\A_J(X^n)$ with $|J| = m$. In particular, $\A_{\age{n}}(X^n)= \A\bigl(h_+(X)^{\otimes n}\bigr)$. This grading is not to be confused with the one given by the dimension of cycles. We have an associated descending filtration $\Fil^\smalldot$ of~$\A(X^n)$, given by
$$
\Fil^r \A(X^n) = \bigoplus_{m=0}^{n-r}\; \A_{\age{m}}(X^n)\, .
$$
This means that the only terms that contribute to $\Fil^r \A(X^n)$ are those coming from submotives $h_{\nu_1}(X) \otimes \cdots \otimes h_{\nu_n}(X)$ involving at least~$r$ factors~$h_0(X)$. Alternatively, a class in~$\A(X^n)$ lies in $\Fil^r \A(X^n)$ if and only if it is a linear combination of classes of the form $\phi_{J,*}(\alpha)$ for subsets $J \subset \{1,\ldots,n\}$ with $n-|J| \geq r$. In particular, if $J \subset \{1,\ldots,n\}$ and $\beta$ is a class in $\Fil^s\A(X^J)$ then $\phi_{J,*}(\beta) \in \Fil^{s+n-|J|}\A(X^n)$.

If $f \colon X \to Y$ is a morphism of smooth connected $k$-varieties and we take $b = f(a)$ as base point on~$Y$, it follows from Proposition~\ref{f_*Prop}(\romno1) that the induced map $(f^n)_* \colon \A(X^n) \to \A(Y^n)$ is a graded map. In particular, it is strictly compatible with the associated filtrations.

\ssection
\ssectlabel{gam0DelFil1}
{\it Remark.\/} --- If $\alpha \in \A(X)$ we have a class $\Delta^{(n)}_{X,*}(\alpha) \in \A(X^n)$. By definition, $\gamma_{X,a}^n(\alpha)$ is the projection of this class onto the summand $\A_{\age{n}}(X^n)$. Hence $\gamma^n_{X,a}(\alpha) = 0$ in~$\A(X^n)$ if and only if $\Delta^{(n)}_{X,*}(\alpha) \in \Fil^1\A(X^n)$.

\ssection
\ssectlabel{delta(n)}
As before, let $X$ be a connected smooth projective $k$-variety with a base point $a \in X(k)$. For $n\geq 1$, consider the morphism $\delta^{(n)} = (\id_{X^{m-1}} \times \Delta_X) \colon X^n \to X^{n+1}$; so $\delta^{(n)}(x_1,\ldots,x_{n-1},x_n) = (x_1,\ldots,x_{n-1},x_n,x_n)$.

If $J \subset \{1,\ldots,n\}$ is a subset with $n \notin J$, the morphism $\delta^{(n)}_* \colon h(X^n) \to h(X^{n+1})$ induces an isomorphism $h_J(X^n) \isomarrow h_J(X^{n+1})$. If $n \in J$, let $\hat{J} = J \cup \{n+1\}$. In this case we have a commutative diagram
$$
\def\normalbaselines{\baselineskip20pt\lineskip3pt\lineskiplimit3pt}
\matrix{
X^{|J|} & \sizedmapright{~\phi_J~}{\sim} & X^J & \sizedmapright{~\phi_J~}{\phi_J} & X^n \cr
\mapdownl{\delta^{(|J|)}} && \mapdownr{\delta^{(n)}|_{X^J}} && \mapdownr{\delta^{(n)}} \cr
X^{|J|+1} & \sizedmapright{~\phi_J~}{\sim} & X^{\hat{J}} & \sizedmapright{~\phi_J~}{\phi_{\hat{J}}} & X^{n+1}\, .
}
$$
It follows that $\delta^{(n)}_* \colon \A(X^n) \to \A(X^{n+1})$ respects the filtrations.

\ssection
\ssectlabel{mtom+1}
{\it Proposition. --- Let $X$ be a connected smooth projective $k$-variety with a base point $a \in X(k)$. Let $n$ be a positive integer.

{\rm (\romno1)} If $\gamma^n_{X,a}(\alpha) = 0$ for some $\alpha \in \A(X)$ then $\gamma^{n+1}_{X,a}(\alpha) = 0$.

{\rm (\romno2)} We have $\Gamma^n(X,a) = 0$ in~$\A(X^n)$ if and only if $\gamma^n_{X,a} \colon \A(X) \to \A(X^n)$ is the zero map.
\par}
\medskip

\Proof
(\romno1)~As remarked in~\ref{gam0DelFil1}, $\gamma^n_{X,a}(\alpha) = 0$ if and only if $\Delta^{(n)}_{X,*}(\alpha) \in \Fil^1 \A(X^n)$. Now use that $\Delta^{(n+1)} = \delta^{(n)} \circ \Delta^{(n)}$ and the fact just explained that $\delta^{(n)}_*$ respects the filtrations.

(\romno2)~Assume that $\Gamma^n(X,a) = 0$ in~$\A(X^n)$ and let $\alpha \in \A(X)$. Because the map~$\gamma^n_{X,a}$ is linear and $\gamma^n_{X,a}[X] = \Gamma^n(X,a)$ by definition, we may assume that $\alpha \in \A_i(X)$ for some $i<\dim(X)$. We know that the class of the small diagonal~$\Delta^{(n)}_X$ lies in~$\Fil^1 \A(X^n)$; this means we can write
$$
\bigl[\Delta^{(n)}_X\bigr] = \sum_{J \subsetneq \{1,\ldots,n\}}\; \beta_J
$$
with $\beta_J \in \A_J(X^n)$. By definition of~$\A_J(X^n)$ we have $\beta_J = \phi_{J,*}(b_J)$ for some class~$b_J$ on~$X^J$. To prove that $\Delta^{(n)}_{X,*}(\alpha) = \bigl[\Delta_X^{(n)}\bigr] \cdot \pr_n^*(\alpha)$ lies in $\Fil^1\A(X^n)$ we now only have to remark that 
$$
\beta_J \cdot \pr_n^*(\alpha) = \phi_{J,*}\Bigl(b_J \cdot (\pr_n \circ \phi_J)^*\bigl(\alpha\bigr)\Bigr)\, ,
$$
and that for $J \subsetneq \{1,\ldots,n\}$ any class in the image of~$\phi_{J,*}$ lies in~$\Fil^1\A(X^n)$. 
\QED
\medskip

For the classes $\Gamma^n(X,a)$ the stability result in~(\romno1) is O'Grady's Proposition~2.4 in~[\ref{OGr}]. As we shall now show, part~(\romno2) of the proposition can be refined. The idea is that we can view $\Gamma^{m+n}(X,a)$ as a correspondence from~$X^m$ to~$X^n$.

\ssection
\ssectlabel{GammaXmToXn}
{\it Proposition. --- Let $X$ be a connected smooth projective $k$-variety with base point $a \in X(k)$. Suppose $m$ and~$n$ are positive integers such that $\Gamma^{m+n}(X,a) = 0$ in $\A(X^{m+n})$. Then
$$
\sum_{\emptyset \neq K \subset \{1,\ldots,m\}}\; (-1)^{|K|} \cdot \gamma_{X,a}^n\Bigl(\Delta_X^{(K),*}(\xi) \Bigr) = 0 \qquad \hbox{in $\A(X^n)$}
$$
for all classes $\xi \in \CH^{>0}(X^m)$. Here $\Delta_X^{(K)} \colon X \to X^m$ denotes the composition of the diagonal $\Delta_X \colon X \to X^K$ and the closed embedding $\phi_K \colon X^K \hookrightarrow X^m$.}
\medskip

\Proof
We may assume $\dim(X) > 0$. By definition,
$$
\Gamma^{m+n}(X,a) = \sum_{\emptyset \neq J \subset \{1,\ldots,m+n\}}\, (-1)^{m+n-|J|} \cdot \bigl[\Delta^{(J)}_X\bigr]\, .
$$
Write the non-empty subsets $J \subset \{1,\ldots,m+n\}$ as $J = K\cup L$ with $K \subset \{1,\ldots,m\}$ and $L = \{m+1,\ldots,m+n\}$. Viewing $\bigl[\Delta_X^{(J)}\bigr]$ as a correspondence from~$X^m$ to~$X^n$, its effect on cycle classes is given by $\xi \mapsto \Delta_{X,*}^{(L)} \bigl(\Delta_X^{(K),*}(\xi) \bigr)$, where in the notation $\Delta_{X,*}^{(L)}$ we treat~$L$ as a subset of $\{1,\ldots,n\}$.

If $K = \emptyset$, the map $\Delta_X^{(K)}$ is the inclusion of the point $(a,\ldots,a)$ in~$X^m$; so $\Delta_X^{(K),*}(\xi) = 0$ for $\xi \in \CH^{>0}(X^m)$. If $K \neq \emptyset$ then
$$
\sum_{L}\, (-1)^{m+n-|K\cup L|}\cdot \Delta^{(K\cup L)}_{X,*}(\xi) = (-1)^{m-|K|} \cdot \gamma_{X,a}^n\Bigl(\Delta_X^{(K),*}(\xi) \Bigr)
$$
and the proposition follows.
\QED

\ssection
\ssectlabel{gammaOnDecompos}
{\it Corollary. --- If $\Gamma^{m+n}(X,a) = 0$ in $\A(X^{m+n})$ then $\gamma^n_{X,a} \colon \A(X) \to \A(X^n)$ is zero on the image of the product map $\A^{>0}(X)^{\otimes m} \to \A(X)$. In particular, if $\Gamma^{n+1}(X,a) = 0$ then $\gamma^n_{X,a}(\xi) = 0$ for all $\xi \in \A^{>0}(X)$.}
\medskip

\Proof
In the proposition, take $\xi = \xi_1 \times \cdots \times \xi_m$ for classes $\xi_i \in \CH^{>0}(X)$. For $K \neq \{1,\ldots,m\}$ we have $\Delta_X^{(K),*}(\xi)=0$. For $K = \{1,\ldots,m\}$ we have $\Delta_X^{(K),*}(\xi) = \xi_1\cdots\xi_m$. Hence we find that $\gamma_{X,a}^n\bigl(\xi_1\cdots\xi_m) = 0$.
\QED

\ssection
\ssectlabel{RelativeDim}
{\it Corollary. --- Let $f \colon X \to Y$ be a surjective morphism of connected smooth projective $k$-varieties. Let $a \in X(k)$ and $b = f(a)$. Let $r = \dim(X) - \dim(Y)$. If $\Gamma^{n+r}(X,a) = 0$ for some $n \geq 1$ then $\Gamma^n(Y,b) = 0$.} 
\medskip

\Proof
There exists a vector bundle~$\cE$ on~$Y$ such that $X$, as a scheme over~$Y$, is isomorphic to a closed subscheme of the projective bundle~$\mP(\cE)$. Let $\ell = c_1\bigl(\cO_{\mP(\cE)}(1)\bigr)$ and write $\ell_X$ for its pull-back to~$X$. We have $f_*(\ell_X^r) \in \CH^0(Y)$; so $f_*(\ell_X^r) = N \cdot [Y]$ for some integer~$N$. By pulling back to the generic point of~$Y$ we see that $N \neq 0$. So by \ref{f_*Prop}(\romno2), $\Gamma^n(Y,b)$ is proportional to $f^{\otimes n}_*\bigl(\gamma^n_{X,a}(\ell_X^r)\bigr)$, which vanishes by Corollary~\ref{gammaOnDecompos}.
\QED

\ssection
\ssectlabel{OGradyConj}
{\it Remark.\/} --- O'Grady has made the conjecture that for a hyperk\"ahler variety~$X$ of dimension~$2n$ there should exist a base point $a \in X(k)$ such that $\Gamma^{2n+1}(X,a) = 0$ in $\CH(X^{2n+1})$. By Corollary~\ref{gammaOnDecompos}, if this is true then the intersection pairing $\Pic(X)^{\otimes 2n} \to \CH_0(X)$ takes values in $\mQ \cdot [a]$. This last property is known (for a suitable choice of $a\in X(k)$) for Hilbert schemes of K3 surfaces and Fano varieties of cubic $4$-folds by results of Voisin~[\ref{VoisHK}], Theorems 1.4(2) and 1.5, and also for generalized Kummer varieties by a result of Fu~[\ref{Fu}], Theorem 1.6.

The vanishing of $\Gamma^n(X,a)$ implies that $n > \dim(X)$; see Theorem~\ref{AlbEstim} below. So for $X$ hyperk\"ahler of dimension~$2n$, the vanishing of $\Gamma^{2n+1}(X,a)$ is the strongest possible result. By Corollary~\ref{RelativeDim}, if O'Grady's conjecture is true then for all varieties~$Y$ dominated by~$X$ we again have the optimal result that $\Gamma^{\dim(Y)+1}(Y,b)$ vanishes in the Chow ring. This suggests that only very special varieties are dominated by a hyperk\"ahler variety. Another indication of this is given by a result of Lin~[\ref{Lin}], Theorem~1.1. He takes for $X$ a Hilbert scheme of points on a complex K3 surface with infinitely many rational curves; then he proves that if $X$ dominates a variety~$Y$ with $\dim(Y) < \dim(X)$, then $Y$ is rationally connected.
\medskip

Part~(\romno2) of the next lemma gives a refinement of the stability result in Proposition~\ref{mtom+1}(\romno1).

\ssection
\ssectlabel{DeltaXFilLem}
{\it Lemma. --- In the situation of\/~{\rm \ref{modaer}}, suppose $\bigl[\Delta^{(n)}_X\bigr] \in \Fil^r \A(X^n)$ for some $r \geq 1$. 

{\rm (\romno1)} For all $i \in \{0,\ldots,r\}$ we have $\bigl[\Delta^{(n-i)}_X\bigr] \in \Fil^{r-i} \A(X^{n-i})$.

{\rm (\romno2)} For all $i \geq 0$ we have $\bigl[\Delta^{(n+i)}_X\bigr] \in \Fil^{r+i} \A(X^{n+i})$.\par}
\medskip

\Proof
In both statements, it suffices to do the case $i=1$. Part~(\romno1) readily follows from the definitions by taking the image of $\bigl[\Delta^{(n)}_X\bigr]$ under a projection $X^n \to X^{n-1}$.

For (\romno2), suppose $\bigl[\Delta^{(n)}_X\bigr] \in \Fil^r \A(X^n)$ with $r \geq 1$. In particular, $\Gamma^n(X,a) = 0$ in~$\A(X^n)$, which by Proposition~\ref{mtom+1}(\romno1) implies that $\Gamma^{n+1}(X,a) = 0$ in~$\A(X^{n+1})$. We can write this as an identity
$$
\bigl[\Delta^{(n+1)}_X\bigr] = \sum_J\; (-1)^{n-|J|} \cdot \bigl[\Delta^{(J)}_X\bigr]
$$
in $\A(X^{n+1})$, where the sum runs over the non-empty subsets $J \subsetneq \{1,\ldots,n+1\}$, and where we recall that $\Delta^{(J)}_X$ is the small diagonal of $X^J$, viewed as a cycle on~$X^{n+1}$. If $|J| \leq n-r$ then it is clear that $\Delta^{(J)}_X \in \Fil^{r+1}\A(X^{n+1})$. If not, then $n+1-r \leq |J| \leq n$ and by the assumption that $\bigl[\Delta^{(n)}_X\bigr] \in \Fil^r \A(X^n)$ together with~(\romno1) the small diagonal on $X^J$ lies in $\Fil^{|J|-n+r}\A(X^J)$. Since $\Delta_X^{(J)}$ is obtained by pushing forward this small diagonal via $\phi_J \colon X^J \hookrightarrow X^{n+1}$, it again follows that $\bigl[\Delta_X^{(J)}\bigr] \in \Fil^{r+1}\A(X^{n+1})$.
\QED
\medskip

We now investigate how changing the base point affects the vanishing of $\Gamma^n(X,a)$.

\ssection
\ssectlabel{IndepBasePt}
{\it Proposition. --- Let $X$ be a connected smooth projective $k$-variety. Let $a$ and~$a^\prime$ be $k$-valued points of~$X$. If $\Gamma^n(X,a) = 0$ in $\A(X^n)$ for some $n > 1$ then $\Gamma^{2n-2}(X,a^\prime) = 0$ in $\A(X^{2n-2})$.}
\medskip

\Proof
Let $\pi^\prime_+ = [\Delta_X] - X\times \{a^\prime\}$ be the projector that cuts out the motive $h_+(X,a^\prime)$. We write it as $\pi^\prime_+ = \pi_+ + X \times \bigl[\{a\}-\{a^\prime\}\bigr]$. This gives
$$
\bigl(\pi^\prime_+\bigr)^{\otimes (2n-2)} = \sum_{J\subset \{1,\ldots,2n-2\}}\; \pi_+^{\otimes J} \,\otimes\, \Bigl(X \times \bigl[\{a\}-\{a^\prime\}\bigr]\Bigr)^{\otimes J^\prime}\, , \eqlabel{pi'-pi^m}
$$
where we write $J^\prime = \{1,\ldots,2n-2\}\setminus J$.

By Corollary~\ref{gammaOnDecompos}, the assumption that $\Gamma^n(X,a) = 0$ implies that $\gamma^m_{X,a}(a^\prime) = 0$ for all $m \geq n-1$. But $\gamma^m_{X,a}(a^\prime) = \bigl[\{a^\prime\}-\{a\}\bigr]^{\otimes m}$; so in \eqref{pi'-pi^m} we may sum only over the subsets $J \subset \{1,\ldots,2n-2\}$ of cardinality $\geq n$. On the other hand, by \ref{DeltaXFilLem}(\romno2) we have $\bigl[\Delta_X^{(2n-2)}\bigr] \in \Fil^{n-1}\A(X^{2n-2})$, which means that $\pi_+^{\otimes J} \otimes \bigl(X \times \bigl[\{a\}-\{a^\prime\}\bigr]\bigr)^{\otimes J^\prime}$ kills $\bigl[\Delta_X^{(2n-2)}\bigr]$ for all index sets~$J$ with $|J| > (2n-2) - (n-1) = n-1$. Together this gives that $\bigl(\pi^\prime_+\bigr)^{\otimes (2n-2)}\bigl[\Delta_X^{(2n-2)}\bigr] = 0$, i.e., $\Gamma^{2n-2}(X,a^\prime) = 0$. 
\QED
\medskip

As an example, on a K3 surface~$X$ with distinguished point class $o_X$ we have $\Gamma^3(X,o_X) = 0$ by [\ref{BeauVoi}], Proposition~3.2. By  Proposition~\ref{IndepBasePt} it follows that for any base point $a \in X(k)$ we have $\Gamma^4(X,a) = 0$, and by Corollary~\ref{gammaOnDecompos} we in fact find that $\Gamma^3(X,a) = 0$ if and only if $a = o_X$ in $\CH_0(X)$.

We finish this section by reproving Proposition~0.2 of O'Grady's paper~[\ref{OGr}], which is an easy consequence of the above.

\ssection
{\it Proposition. ---  Let $X$ and~$Y$ be connected smooth projective $k$-varieties with base points $a \in X(k)$ and $b \in Y(k)$. Suppose that $\Gamma^m(X,a) = 0$ in~$\A(X^m)$ and $\Gamma^n(Y,b) = 0$ in~$\A(Y^n)$ for some positive integers $m$ and~$n$. Then $\Gamma^{m+n-1}\bigl(X\times Y,(a,b)\bigr) = 0$ in $\A\bigl((X\times Y)^{m+n-1}\bigr)$.\par}
\medskip

\Proof
By Lemma~\ref{DeltaXFilLem}(\romno2) we have 
$$
\bigl[\Delta_X^{(m+n-1)}\bigr] \in \Fil^n\A(X^{m+n-1})\, ,\qquad \bigl[\Delta_Y^{(m+n-1)}\bigr] \in \Fil^m\A(X^{m+n-1})\, .
$$ 
This means we can write $\bigl[\Delta_X^{(m+n-1)}\bigr] = \sum_J\, \phi_{X,J,*}(\alpha_J)$, where the sum runs over the subsets $J \subset \{1,\ldots,m+n-1\}$ of cardinality at most $m-1$, and where $\alpha_J$ is a class on~$X^J$. Similarly, $\bigl[\Delta_Y^{(m+n-1)}\bigr] = \sum_K\, \phi_{Y,K,*}(\beta_K)$, where the subsets $K \subset \{1,\ldots,m+n-1\}$ have cardinality at most $n-1$ and $\beta_K \in \A(Y^K)$.

Writing $p \colon (X\times Y)^{m+n-1} \to X^{m+n-1}$ and $q\colon (X\times Y)^{m+n-1} \to Y^{m+n-1}$ for the projections, 
$$
\bigl[\Delta_{X\times Y}^{m+n-1}\bigr] = p^*\bigl[\Delta_X^{(m+n-1)}\bigr] \cdot q^*\bigl[\Delta_Y^{(m+n-1)}\bigr] \, .
$$
Given subsets $J$, $K \subset \{1,\ldots,m+n-1\}$ of cardinality at most $m-1$ and $n-1$, respectively, there is an index $\nu \in \{1,\ldots,m+n-1\}$ that is not in $J \cup K$. Setting $L = \{1,\ldots,\hat{\nu},\ldots,m+n-1\}$ it is then clear that $p^*\phi_{X,J,*}(\alpha_J) \cdot q^*\phi_{Y,K,*}(\beta_K)$ is a class in the image of the push-forward under $\phi_{X\times Y,L} \colon (X\times Y)^L \longhookrightarrow (X\times Y)^{m+n-1}$. Hence, $\bigl[\Delta_{X\times Y}^{m+n-1}\bigr] \in \Fil^1\A\bigl((X\times Y)^{m+n-1}\bigr)$, which means that $\Gamma^{m+n-1}(X\times Y) = 0$.
\QED

\section{Vanishing results on curves and abelian varieties}{CAV}
\bigskip

\noindent
We begin by recalling a result of O'Grady~[\ref{OGr}] about the vanishing of modified diagonals in cohomology.

\ssection
\ssectlabel{AlbEstim}
{\it Theorem. --- Let $X$ be a connected smooth projective $k$-variety with base point $a \in X(k)$. Let $d = \dim(X)$ and let $e$ be the dimension of the image of the Albanese map $\alb \colon X \to \Alb(X)$. Then $\Gamma^n(X,a) \sim_\hom 0$ if and only if $n > d+e$.
\par}
\medskip

\Proof
Let $H^\smalldot$ denote $\ell$-adic cohomology for some prime $\ell \neq \char(k)$. Throughout, we view $H^\smalldot(X) = \oplus_{i=0}^{2d}\, H^i(X)$ as a superspace; in particular, if $i$ is odd then $\Sym^m(H^i)$ has $\wedge^m H^i(X)$ as its underlying vector space.

The cohomology class $[\Gamma^n]$ of $\Gamma^n(X,a)$ lies in the degree $2d(n-1)$-part of $\Sym^n\bigl(H^\smalldot(X)\bigr)$. We have
$$
\Sym^n\bigl(H^\smalldot(X)\bigr) = \bigoplus_{m=(m_0,\ldots,m_d)\atop |m|=n}\; \bigotimes_{j=0}^{2d}\, \Sym^{m_j}\bigl(H^j(X)\bigr)\, ,
$$
where the summand $S(m) = \otimes_j\, \Sym^{m_j}\bigl(H^j(X)\bigr)$ lies in degree $\sum_{j=0}^{2d}\, j \cdot m_j$. By Remark~\ref{h0h*h2d} we know that the component of $[\Gamma^n]$ in $S(m)$ is zero if $m_0 > 0$ or $m_{2d}>0$. Next consider a sequence $m = (m_0,m_1,\ldots,m_{2d})$ with $m_0 = m_{2d} = 0$. The component of~$[\Gamma^n]$ in~$S(m)$ is then the same as the component of the cohomology class of the small diagonal $\Delta_X^{(n)}$. If $\mu = (m_{2d},m_{2d-1},\ldots,m_0)$ is the reverse sequence, the intersection pairing on $\Sym^n\bigl(H^\smalldot(X)\bigr) \subset H^\smalldot(X^n)$ restricts to a perfect pairing $S(m) \times S(\mu) \to k$, and for $m^\prime \neq \mu$ the pairing $S(m) \times S(m^\prime) \to k$ is zero. For $\beta \in S(\mu)$ we have $\bigl[\Delta_X^{(n)}\bigr] \cdot \beta = \deg\bigl(\Delta^*(\beta)\bigr)$, and we claim that this is zero whenever $m_{2d-1} > 2e$. Assuming this for a moment, the ``if'' statement in the theorem follows, as the highest degree we can get under the restrictions $m_0=m_{2d}=0$ and $m_{2d-1} \leq 2e$ is $2e(2d-1) + (n-2e)(2d-2) = 2e+2nd-2n$, so that for $n > d+e$ we cannot reach degree $2d(n-1)$.

It remains to be shown that for $i > 2e$ the multiplication map $\Delta^* \colon \Sym^i H^1(X) \to H^i(X)$ is zero. For this we use that $H^1(\alb) = \alb^* \colon H^1(\Alb_X) \to H^1(X)$ is an isomorphism. We have a commutative diagram 
$$
\epsfbox{Diagram.1}
$$ 
But $H^i(\alb)$ factors through $H^i\bigl(\alb(X)\bigr)$, which is zero for $i>2e$. 

Finally we show that $\Gamma^{d+e}(X,a)$ is not homologically trivial, which by Proposition~\ref{mtom+1}(\romno1) gives the ``only if'' in the theorem. The only sequence $m = (m_0,m_1,\ldots,m_{2d})$ with $|m|=d+e$ and $m_0 = m_{2d} = 0$ that reaches degree $2d(d+e-1)$ is $m = (0,\ldots,0,d-e,2e,0)$. With $\mu$ the reverse sequence, it suffices to produce an element $\beta \in S(\mu) = \Sym^{2e} H^1(X) \otimes \Sym^{d-e} H^2(X)$ for which $\Delta^*(\beta)$ has degree $\neq 0$. For this we take polarizations $L_1 \in H^2(\Alb_X)$ and $L_2 \in H^2(X)$; then take $\beta = \Sym^{2e}H^1(\alb)\bigl(L_1^e\bigr) \otimes L_2^{d-e}$. Because the map $H^{2e}(\alb)$ is injective, $\Delta^*(\beta)$ has positive degree and we are done. 
\QED
\medskip

Next we turn to abelian varieties. The result we prove was conjectured by O'Grady in the first version of~[\ref{OGr}]. He also proved it for $g\leq 2$.

\ssection
\ssectlabel{AVResult}
{\it Theorem. --- Let $X$ be an abelian variety of dimension~$g$ over a field~$k$. Let $a \in X(k)$ be a base point. Then $\Gamma^n(X,a)=0$ in $\CH(X^n)$ for all $n > 2g$.
\par}
\medskip

\Proof
We give $X$ the group structure for which~$a$ is the origin. For $m \in \mZ$ let $\mult(m) \colon X \to X$ be the endomorphism given by multiplication by~$m$. By [\ref{DenMur}], Corollary~3.2, we have a motivic decomposition $h(X) = \oplus_{i=0}^{2g}\,  h_i(X)$ in~$\Mot_k$ that is stable under all endomorphisms~$\mult(m)_*$, and such that $\mult(m)_*$ is multiplication by~$m^i$ on~$h_i(X)$. (The result is stated in op.\ cit.\ for the cohomological theory but is easily transcribed into the homological language.) The relation with~\eqref{eq:h0h+} is that $h_0(X,a) = h_0(X)$ and $h_+(X,a) = \bigoplus_{i>0}\, h_i(X)$.

For $n \geq 1$ this induces a decomposition 
$$
h(X^n) = \bigoplus_{\bi = (i_1,\ldots,i_n)}\quad \bigotimes_{j=1}^n\; h_{i_j}(X) \, ,  
$$ 
where the sum runs over the elements $\bi = (i_1,\ldots,i_n)$ in $\{0,\ldots,2g\}^n$. Under this decomposition we have
$$
h_r(X^n) = \bigoplus_{|\bi|=r}\quad \bigotimes_{j=1}^n\; h_{i_j}(X) \, , 
$$
where the sum runs over the $n$-tuples $\bi$ with $|\bi| = i_1+\cdots+i_n$ equal to~$r$.

Now observe that $\bigl[\Delta^{(n)}_X\bigr] \in \CH\bigl(h_{2g}(X^n)\bigr)$, because $\mult(m)_* \bigl[\Delta^{(n)}_X\bigr] = m^{2g} \cdot \bigl[\Delta^{(n)}_X\bigr]$ for all~$m$. The theorem follows, since for $n>2g$ and $\bi = (i_1,\ldots,i_n)$ in $\{0,\ldots,2g\}^n$ with $|\bi| = 2g$ there is at least one index~$j$ with $i_j = 0$. 
\QED
\medskip

Next we turn to curves. Part~(\romno1) of the next result is due to Gross and Schoen; see [\ref{GrSch}], Proposition~3.1. This result is also an immediate consequence of Theorem~\ref{AlbEstim}. Part (\romno2) is due to Polishchuk; see~[\ref{Polish}], Corollary~4.4(\romno4). Part~(\romno3) is essentially due to Polishchuk and the first author in~[\ref{BMAP}] (see especially the proof of loc.\ cit.\ Theorem~8.5) but we need to combine the calculations that are done there with some known facts about the Chow ring of the Jacobian, as we shall now explain.

\ssection
\ssectlabel{CurveResult}
{\it Theorem. --- Let $C$ be a complete nonsingular curve of genus~$g$ over a field~$k$ with a base point $a \in X(k)$. Then
\item{\rm (\romno1)} $\Gamma^n(C,a) \sim_\hom 0$ for all $n > 2$;
\item{\rm (\romno2)} $\Gamma^n(C,a) \sim_\alg 0$ (modulo torsion) for all $n > \gon(C)$;
\item{\rm (\romno3)} $\Gamma^n(C,a) = 0$ in $\CH_1(C^n)$ for all $n > g+1$.
\par}
\medskip

\Proof
For curves of genus~$0$ the result is trivial. (Because we work modulo torsion, we may extend the ground field and assume $C=\mP^1$; then note that the diagonal of $\mP^1 \times \mP^1$ is rationally equivalent to $\bigl(\{\pt\}\times \mP^1\bigr) + \bigl(\mP^1 \times \{\pt\}\bigr)$.) Hence we may assume $g>0$. Let $\iota \colon C \to J$ be the closed embedding associated with the base point~$a$. As discussed above, $h(J) = \oplus_{i=0}^{2g}\, h_i(J)$. This means we can decompose $\bigl[\iota(C)\bigr] \in \CH_1(J)$ as
$$
\bigl[\iota(C)\bigr] = \sum_{i=0}^{2g}\, \gamma_i
$$
with $\gamma_i \in \CH\bigl(h_i(J)\bigr)$. In particular, for $m\in \mZ$ we have $\mult(m)_*(\gamma_i) = m^i \cdot \gamma_i$. It is known that: 
\itemitem{(a)} $\gamma_i \neq 0$ only for $i \in \{2,\ldots,g+1\}$;
\itemitem{(b)} $\gamma_i$ is torsion modulo algebraic equivalence for $i> \gon(C)$;
\itemitem{(c)} $\gamma_i$ is homologically trivial for $i \neq 2$. 

\noindent
In fact, (c) holds because $\mult(m)_*$ acts on $H^{2g-2}(J)$ as multiplication by~$m^2$, (a) follows from the precise summation range in the main theorem of~[\ref{BeauvChow}] (in the notation of~[\ref{BeauvChow}] our~$\gamma_i$ lies in $\CH^{g-1}_{i-2}(J)$), and (b) is a result of Colombo and van Geemen~[\ref{ColvG}].

We denote by $C^{[d]}$ the $d$th symmetric power of~$C$ and let $C^{[\bullet]} = \coprod_{d\geq 0}\, C^{[d]}$, which is a monoid scheme. Let $\CH\bigl(C^{[\bullet]}\bigr) = \oplus_{d \geq 0}\, \CH\bigl(C^{[d]}\bigr)$, which is a $\mQ$-algebra for the Pontryagin product. The maps $u_d\colon C^{[d]} \to J$ give us a morphism $u \colon C^{[\bullet]} \to J$, which induces a homomorphism $u_* \colon \CH\bigl(C^{[\bullet]}\bigr)\to \CH(J)$. By [\ref{BMAP}], Theorem~3.4, there is a $\mQ$-subalgebra $\mK \subset \CH\bigl(C^{[\bullet]}\bigr)$ such that the restriction of~$u_*$ to~$\mK$ gives an isomorphism $\mK \isomarrow \CH(J)$. Further, by ibid., Lemma~8.4 and the proof of Theorem~8.5, all classes $\Gamma^n(C,a)$ lie in this subalgebra~$\mK$ and we have, for $n\geq 2$,
$$
u_*\bigl(\Gamma^n(C,a)\bigr) = n! \cdot \sum_{i=0}^{2g}\; S(i,n) \cdot \gamma_i\, , 
$$
where $S(i,n)$ denotes the Stirling number of the second kind. Note that $S(i,n)=0$ if $n>i$. Putting together these facts, the theorem follows from (a)--(c) above.
\QED

\ssection
Let us now discuss to what extent the above results are sharp. 

For abelian varieties, our result in Theorem~\ref{AVResult} is sharp, since by Theorem~\ref{AlbEstim} $\Gamma^{2g}(X,a)$ is not even homologically trivial. The same remark applies to part~(\romno1) of Theorem~\ref{CurveResult}.

Part~(\romno2) of Theorem~\ref{CurveResult} is conjecturally sharp for the generic curve $C$ of genus~$g$. In fact, under the genericity assumption it is expected that $\gamma_i$ is not algebraically trivial for $i = \lfloor (g + 3)/2 \rfloor = \gon(C)$. We refer to [\ref{VoisInfi}] for recent results (in characteristic~$0$) towards this conjecture.

Finally, (\romno3) of Theorem~\ref{CurveResult} is sharp for the generic pointed curve in characteristic~$0$. This is proven by the second author in~[\ref{QYThesis}], Proposition~5.14, which gives $\gamma_{g+1} \neq 0$.

\section{Double covers}{Double}
\bigskip

\noindent
The following result proves a conjecture made by O'Grady in~[\ref{OGr}]. We had originally hoped to extend this to more general covers, but our method leads to some non-trivial combinatorial problems. As we just learned that Claire Voisin has obtained such a more general result using a different argument, we restrict ourselves to double covers. As in Section~\ref{Firstvanish}, we consider an adequate equivalence relation~$\sim$ and write $\A(X) = \CH(X)/{\sim}$.

\ssection
\ssectlabel{DoubleThm}
{\it Theorem. --- Let $f\colon X \to Y$ be a double cover. Let $\sigma$ be the corresponding involution of~$X$. Let $a \in X(k)$ be a base point such that $a \sim \sigma(a)$, and write $b = f(a)$. If $\Gamma^n(Y,b) = 0$ in~$\A(Y^n)$ then $\Gamma^{2n-1}(X,a) = 0$ in~$\A(X^{2n-1})$.
\par}
\medskip

\ssection
As a preparation for the proof we need to introduce some notation. Given an integer~$m$ and a subset $J \subset \{1,\ldots,m\}$, let $Z_J \subset X^m$ denote the image of the morphism $\zeta_J \colon X \to X^m$ for which
$$
\pr_j \circ \zeta_J = \cases{\sigma & if $j \in J$;\cr \id_X & if $j \notin J$.}
$$
If $J^\prime$ is the complement of~$J$, we have $Z_{J^\prime} = Z_J$. Further, $Z_\emptyset = Z_{\{1,\ldots,m\}} = \Delta_X^{(m)}$. 

For $r \leq m$, let
$$
V_r = \sum_{{J \subset \{1,\ldots,m\} \atop |J|=r}}\; \bigl[Z_J\bigr]\, .
$$
It follows from the previous remarks that $V_{m-r}=V_r$ and that $V_0 = V_m = \bigl[\Delta_X^{(m)}\bigr]$. We write $V_r^{(m)}$ if there is a need to specify~$m$.

The pull-back of the class $\bigl[\Delta_Y^{(m)}\bigr]$ is ${1\over 2} \cdot \sum_{r=0}^m\, V_r$.

\ssection
For $(i,j) \in \{1,\ldots,m\} \times \{1,\ldots,m+1\}$, consider the morphism $\phi_{i,j} \colon X^m \to X^{m+1}$ given by
$$
(x_1,\ldots,x_m) \mapsto \bigl(x_1,\ldots,x_{j-1},\sigma(x_i),x_j,\ldots,x_m\bigr)\, .
$$
Let $\Phi$ be the sum of the graphs of the~$\phi_{i,j}$; so, $\Phi = \sum_{i,j}\, \bigl[\Gamma_{\phi_{i,j}}\bigr]$. This is a correspondence of degree~$0$ from~$X^m$ to~$X^{m+1}$. Again we write $\Phi^{(m)}$ if we want to specify~$m$.

\ssection
\ssectlabel{Phi*V}
{\it Lemma. --- For $r \leq m$ we have
$$
\Phi_*\bigl(V_r^{(m)}\bigr) = r(m+1-r) \cdot V_r^{(m+1)} + (r+1)(m-r) \cdot V_{r+1}^{(m+1)}\, .
$$
\par}
\vskip-\lastskip\medskip

\Proof
Given $j\in \{1,\ldots,m\}$, let $\alpha_j \colon \{1,\ldots,m\} \to \{1,\ldots,m+1\}$ be the strictly increasing map such that $j$ is not in the image of~$\alpha_j$. Fix some subset $K \subset \{1,\ldots,m+1\}$. We have to count the number of choices for $J \subset \{1,\ldots,m\}$ with $|J|=r$ and an index pair $(i,j)$ as above such that $\phi_{i,j,*}[Z_J] = [Z_K]$. It is clear that there are no such choices unless $|K|=r$ or $|K|=r+1$. If $|K|=r$ then we can choose $j \notin K$ and $i \in \alpha_j^{-1}(K)$ arbitrarily; once these choices are made there is a unique $J \subset \{1,\ldots,m\}$ with $|J|=m$ such that $\phi_{i,j,*}[Z_J] = [Z_K]$. Note that the number of choices in this case is $(m+1-r)r$. Similarly, if $|K|=r+1$ we have to choose $j \in K$ and $i \notin \alpha_j^{-1}(K)$ and then there is again a unique choice for~$J$ such that $\phi_{i,j,*}[Z_J] = [Z_K]$. In this case the number of choices is $(r+1)(m-r)$. 
\QED

\ssection
\ssectlabel{LemmaVr}
{\it Lemma. --- Notation and assumptions as in Theorem\/~{\rm \ref{DoubleThm}}. If $\Gamma^n(Y,b) = 0$ in $\A(Y^n)$ then $\sum_{r=0}^{m+n}\, r^j(m+n-r)^j \cdot V_r^{(m+n)}$ lies in $\Fil^1\A(X^{m+n})$ for all $m\geq j \geq 0$.
\par}
\medskip

\Proof
We use induction on~$m$. For $m=0$ the assumption that $\Gamma^n(Y,b) = 0$ means that $\bigl[\Delta_Y^{(n)}\bigr] \in \Fil^1\A(Y^n)$. Pulling back to~$X^n$ and using that $a \sim \sigma(a)$ we find that $\sum_{r=0}^n\, V_r$ lies in~$\Fil^1\A(X^n)$.

Assuming the assertion is true for some~$m$, let us prove it for $m+1$. By Proposition~\ref{mtom+1}, $\bigl[\Delta_Y^{(n)}\bigr] \in \Fil^1\A(Y^n)$ implies that $\bigl[\Delta_Y^{(n+1)}\bigr] \in \Fil^1\A(Y^{n+1})$. So the assertion for $j < m+1$ follows from the induction hypothesis, replacing $n$ with $n+1$.

It remains to consider the case $j=m+1$. Let
$$
W = \Phi^{(m+n)}_*\Biggl(\sum_{r=0}^{m+n}\, r^m(m+n-r)^m \cdot V_r^{(m+n)}\Biggr)\, .
$$
By the induction assumption, $\sum_{r=0}^{m+n}\, r^m(m+n-r)^m \cdot V_r^{(m+n)}$ lies in $\Fil^1\A(X^{m+n})$, and by the same argument as in~\ref{delta(n)}, $\Phi_* \colon \A(X^n) \to \A(X^{n+1})$ respects the filtrations; hence, $W \in \Fil^1\A(X^{m+n+1})$.

By Lemma~\ref{Phi*V}, $W$ equals
$$
\eqalign{
&\sum_{r=0}^{m+n}\; r^m (m+n-r)^m \cdot \Bigl(r(m+n+1-r) \cdot V_r^{(m+n+1)} + (r+1)(m+n-r)\cdot V_{r+1}^{(m+n+1)}\Bigr)\cr
&= \sum_{s=0}^{m+n+1}\; \Bigl(s^{m+1}(m+n-s)^m(m+n+1-s) + (s-1)^m s (m+n+1-s)^{m+1}\Bigr) \cdot V_s^{(m+n+1)}\cr
&= \sum_{s=0}^{m+n+1}\; s(m+n+1-s) \cdot \Bigl(s^m(m+n-s)^m + (s-1)^m(m+n+1-s)^m\Bigr)  \cdot V_s^{(m+n+1)}\, .\cr
}
$$
Putting $x = s$ and $y = m+n+1-s$ we have 
$$
\eqalign{s^m(m+n-s)^m + (s-1)^m(m+n+1-s)^m &= x^m(y-1)^m + (x-1)^my^m \cr
&= \sum_{j=0}^m\; (-1)^j {m\choose j}\cdot (x^j+y^j) x^{m-j}y^{m-j}\, .}
$$
As $x+y = m+n+1$ is constant, we can rewrite this as $2x^my^m + \sum_{j=0}^{m-1}\, c_j \cdot x^jy^j$ for some constants $c_0,\ldots,c_{m-1}$. Hence
$$
\eqalign{
W = {} &2\cdot \sum_{s=0}^{m+n+1}\; s^{m+1}(m+n+1-s)^{m+1} \cdot V_s^{(m+n+1)}\cr
&\qquad + \sum_{j=1}^m c_{j-1} \Biggl( \sum_{s=0}^{m+n+1}\; s^j(m+n+1-s)^j \cdot V_s^{(m+n+1)}\Biggr)\cr}\, .
$$
As we have already shown that for $j<m+1$ 
$$
\sum_{s=0}^{m+n+1}\; s^j(m+n+1-s)^j \cdot V_s^{(m+n+1)} \;\in \Fil^1\A(X^{m+n+1})\, ,
$$
the same is true for the remaining term, i.e., for $j=m+1$.
\QED

\medskip

\Pf{Proof of Theorem\/ {\rm \ref{DoubleThm}}}
Taking $m=n-1$ in Lemma~\ref{LemmaVr} and using that $V_r^{(2n-1)} = V_{2n-1-r}^{(2n-1)}$, we find that
$$
\sum_{r=0}^{n-1}\; \bigl(r(2n-1-r)\bigr)^j \cdot V_r^{(2n-1)} = 0 \quad \hbox{in $\A(X^{2n-1})/\Fil^1$}
$$
for all $j \in \{0,1,\ldots,n-1\}$. The $n \times n$ matrix 
$$
\Bigl(\bigl(r(2n-1-r)\bigr)^j\Bigr)_{r,j=0,\ldots,n-1}
$$
is a Vandermonde matrix with distinct entries in the second column $(j=1)$. Therefore, $\bigl[\Delta_X^{(2n-1)}\bigr] = V_0^{(2n-1)} \in \Fil^1\A(X^{2n-1})$, which means that $\Gamma^{2n-1}(X,a) = 0$.
\QED
\medskip

We refer to the forthcoming paper~[\ref{Voisin}] of Claire Voisin for a generalization of this result to covers of arbitrary degree.

\vskip2.0\bigskipamount plus 2pt minus 1pt%
\goodbreak

\noindent
{\bf References}
\nobreak\vskip.75\bigskipamount plus 2pt minus 1pt%

{\eightpoint

\bibitem{YAIntro}
Y.~Andr\'e, {\it Une introduction aux motifs.\/} Panoramas et synth\`eses 17, Soc.\ Math.\ France 2004. 

\bibitem{BeauvChow}
A.~Beauville, {\it Sur l'anneau de Chow d'une vari\'et\'e ab\'elienne.\/}
Math.\ Ann.\ 273 (1986), 647--651.

\bibitem{BeauVoi}
A.~Beauville, C.~Voisin, {\it On the Chow ring of a K3 surface.\/} J.\ Alg.\ Geom.\ 13 (2004), 417--426.

\bibitem{ColvG}
E.~Colombo, B.~van Geemen, {\it Note on curves in a Jacobian.\/} Compositio Math.\ 88 (1993), 333--353. 

\bibitem{DenMur}
C.~Deninger, J.~Murre, {\it Motivic decomposition of abelian schemes and the Fourier transform.\/} J.\ reine angew.\ Math.\ 422 (1991), 201--219. 

\bibitem{Fu}
L.~Fu, {\it Beauville-Voisin conjecture for generalized Kummer varieties.\/} International Math.\ Res.\ Notices (2014). DOI:10.1093/imrn/rnu053.

\bibitem{GrSch}
B.~Gross, C.~Schoen, {\it The modified diagonal cycle on the triple product of a pointed curve.\/} Ann.\ Inst.\ Fourier 45 (1995), 649--679.

\bibitem{Lin}
H.-Y.~Lin, {\it Rational maps form punctual Hilbert schemes of K3 surfaces.\/} arXiv:1311.0743.

\bibitem{BMAP}
B.~Moonen, A.~Polishchuk, {\it Algebraic cycles on the relative symmetric powers and on the relative Jacobian of a family of curves. II.} J.\ Inst.\ Math.\ Jussieu 9 (2010) 799--846.

\bibitem{MY}
B.~Moonen, Q.~Yin, {\it On a question of O'Grady about modified diagonals.\/} arXiv:1311.1185. 

\bibitem{Polish}
A.~Polishchuk, {\it Algebraic cycles on the relative symmetric powers and on the relative Jacobian of a family of curves. I.} Selecta Math.\ (N.S.) 13 (2007), 531--569. 

\bibitem{OGr}
K.~O'Grady, {\it Computations with modified diagonals.\/} arXiv:1311.0757.

\bibitem{VoisHK}
C.~Voisin, {\it On the Chow ring of certain algebraic hyper-K\"ahler manifolds.\/} Pure Appl.\ Math.\ Q.\ 4 (2008), 613--649. 

\bibitem{VoisInfi}
C.~Voisin, {\it Infinitesimal invariants for cycles modulo algebraic equivalence and $1$-cycles on Jacobians.\/} Alg.\ Geom.\ 2 (2014), 140--165.

\bibitem{Voisin}
C.~Voisin, {\it Some new results on modified diagonals.\/} In preparation.

\bibitem{QYThesis}
Q.~Yin, {\it Tautological cycles on curves and Jacobians.\/} Thesis, Radboud University Nijmegen, 2014. Available at http://www.math.ethz.ch/$\sim$yinqi.

\vskip 2\bigskipamount

\noindent
Radboud University Nijmegen, IMAPP, PO Box 9010, 6500GL Nijmegen, The Netherlands. 

\noindent
b.moonen@science.ru.nl
\medskip

\noindent
ETH Z\"urich, Departement Mathematik, 8092 Z\"urich, Switzerland. 

\noindent
qizheng.yin@math.ethz.ch
\par}

\bye